\documentclass[reqno]{amsart}
\usepackage{amssymb}
\usepackage[matrix,arrow,tips,curve,ps]{xy}
\SelectTips{cm}{10}
\usepackage{enumerate}
\newtheorem{theorem}{Theorem}
\newtheorem{lemma}[theorem]{Lemma}
\newtheorem{proposition}[theorem]{Proposition}
\newtheorem{corollary}[theorem]{Corollary}

\theoremstyle{definition}
\newtheorem{remark}[theorem]{Remark}

\newcommand{\PP}{\mathbb{P}}
\newcommand{\C}{\mathbb{C}}
\renewcommand{\le}{\leqslant}
\renewcommand{\ge}{\geqslant}
\newcommand{\I}{\mathcal{I}}
\newcommand{\M}{\mathcal{M}}
\newcommand{\A}{\mathcal{A}}
\newcommand{\F}{\mathcal{F}}
\newcommand{\cL}{\mathcal{L}}
\newcommand{\FF}{\mathbb{F}}
\newcommand{\cO}{\mathcal{O}}

\newlength{\spazio}

\DeclareMathOperator{\Sym}{Sym}

\begin{document}
\title[Rank of the 2nd Gaussian map for general curves]%
{The rank of the 2nd Gaussian map\\ for general curves}

\author{Alberto Calabri, Ciro Ciliberto \and Rick Miranda}

\email{calabri@dmsa.unipd.it}
\curraddr{DMMMSA, Universit\`a di Padova\\
Via Trieste 63, 35121 Padova, Italy}

\email{cilibert@mat.uniroma2.it}
\curraddr{Dipartimento di Matematica,
Universit\`a di Roma Tor Vergata,
Via della Ricerca Scientifica, 00133 Roma, Italy}

\email{Rick.Miranda@ColoState.edu}
\curraddr{Department of Mathematics,
101 Weber Building,
Colorado State University,
Fort Collins, CO 80523-1874, USA}

\thanks{
\\
The first two authors are partially supported
by G.N.S.A.G.A.-I.N.d.A.M.\ ``Francesco Severi''.}

\begin{abstract}
We prove that, for the general curve of genus $g$, the $2$nd Gaussian map $\mu$
is injective if $g \le 17$ and surjective if $g \ge 18$.
The proof relies on the study of the limit of $\mu$ 
when the general curve of genus $g$ degenerates to a general stable binary curve,
i.e.\ the union of two rational curves meeting at $g+1$ points.
\end{abstract}

\maketitle

\section*{Introduction}
Let $X$ be a smooth, projective curve of genus $g$ 
and let $\cL$ be a line bundle on $X$. Consider the product $X\times X$,  with the projections $p_1,p_2$ to the factors, and the natural morphism $p$ to the symmetric product $X(2)$. One has
$p_* (p_1^* \cL\otimes p_2^*\cL)=\cL^+\oplus \cL^-$, where $\cL^\pm$ are the invariant and anti-invariant line bundles with respect to the involution $(x,y)\mapsto (y,x)$. One has $H^0(\cL^+)\cong {\Sym}^2H^0 (\cL)$ and $H^0 (\cL^-)\cong \wedge ^2H^0 (\cL)$. Restriction to the diagonal of $X(2)$ gives rise to two maps 
\[
\mu_{\cL,1}: {\Sym}^2H^0 (\cL)\to H^0(\cL^{\otimes 2}), \quad w_{\cL,1}: \wedge^2H^0 (\cL)\to H^0(\cL^{\otimes 2}\otimes K_X), 
 \]
 where $K_X$ is the canonical bundle of $X$. Both maps have a well known geometric meaning. The  former is given by considering the map $\phi_\cL: X \to \PP^r:=\PP(H^0(\mathcal L))^*$ defined by the complete linear series determined by $\cL$ and by pulling back to $X$ forms of degree two in $\PP^r$. 
The latter is  given by considering the  composition $\gamma$ of $\phi_\cL$ with the \emph{Gauss map} of $X$ to the  Grassmannian of lines $\mathbb G(1,r)$ and by pulling back to $X$ via $\gamma$ forms of degree one in
$\PP^{\binom{r+1}{2}-1}$.  

The maps $\mu_{\cL,1}$ and $w_{\cL,1}$ are the first instances of two hierarchies of maps $\mu_{\cL,k}$ and $w_{\cL,k}$, defined for all positive integers $k$, and called by some authors  \emph{higher Gaussian maps} of $X$. They are inductively defined by iterated restrictions  to the diagonal of $X(2)$.
Precisely for all $k\ge 2$ one has
 \[
\mu_{\cL,k}:  \ker(\mu_{\cL,k-1}) \to H^0(\cL^{\otimes 2} \otimes K_X^{\otimes 2(k-1)}), \quad w_{\cL,k}: \ker(w_{\cL,k-1}) \to H^0(\cL^{\otimes 2} \otimes K_X^{\otimes (2k-1)}).
 \]
 
These maps are particularly interesting when
$\cL\cong K_X$, in which case we will simply denote them as  $\mu_{k}$ and $w_{k}$. They are both defined at a general point of the moduli space of curves $\M_g$ and it is natural to
guess that they have some modular meaning. 
Indeed,  $\mu_1$ is the codifferential, at the point corresponding to $X$,  of the Torelli map $\tau: \M_g\to \mathcal A_g$, and  Noether's theorem says it is surjective if and only if  $X$ is  non-hyperelliptic. 
 
 The map $w_1$ is called the \emph{Wahl map}, and it is related to important deformation and extendability properties of the canonical image of the curve (cf.\ \cite {BM, W}).
Because of this, it has been studied by various authors, too many to be quoted here.
One the most interesting results concerning it is perhaps a theorem
first proved by Ciliberto, Harris and Miranda in  \cite{CHM}, to the effect
that $w_1$ is surjective, as expected, for a general curve of genus $g=10$ and $g\ge12$.
Moreover, this map is injective, as expected, for a general curve of genus $g\le 8$, cf.\ \cite{CM2}.
Unexpectedly, the Wahl map is not of maximal rank for a general curve of genus $g=9,11$.

In general, all maps $\mu_{k}$ and $w_{k}$
are supposed to be meaningful in the geometry of curves,
especially of curves with general moduli. 
Here we will look in particular at the map $\mu_2\colon \I_2(K_X) \to H^0(X, K_X^{\otimes4})$, 
where $\I_2(K_X)$ is the vector space of forms of degree two
vanishing on the canonical model of $X$.
From now on we will simply denote this map by $\mu$,  
and we will call it the
\emph{2nd Gaussian map} of $X$. 
This map was first considered by Green-Griffiths  in ~\cite{G}
and its importance resides in the fact that it is related to the
2nd fundamental form of the moduli space of curves $\M_g$,
embedded in $\A_g$ via the Torelli map, cf.\ \cite{CPT, CF1, CF2}. 

Despite the unexpected behaviour of the Wahl map for genus $g=9,11$,
a reasonable working hypothesis is that the  2nd Gaussian map $\mu$ should be of maximal rank
for a general curve of any genus $g$.
A dimension count shows that this is equivalent to say that $\mu$ should be injective for a
general curve of genus $g \le 17$ and surjective if $g \ge 18$.
So far, the best result in this direction has been proved by Colombo and Frediani in \cite{CF3},
where, by studying hyperplane sections of high genera of K3 surfaces,
they show that $\mu$ is surjective for a general curve of genus $g>152$.
For other interesting results concerning $\mu$, see also \cite{CF2, CFP}.

In this paper, we prove the maximal rank property for every genus:

\begin{theorem}\label{thm1}
The 2nd Gaussian map $\mu\colon \I_2(K_X) \to H^0(X, K_X^{\otimes4})$
for $X$ a general curve of any genus $g$ has maximal rank,
namely it is injective for $g \le 17$ and surjective for $g \ge 18$.
\end{theorem}

As shown in \cite {CPT}, the map $\mu$ has a lifting
$\rho: \mathcal I_2(K_X)\to {\rm Sym}^2(H^0(K_X^{\otimes 2}))$, which is the
datum of the second fundamental form of the Torelli embedding at the
point corresponding to $X$ in the non-hyperelliptic case. As proved in \cite {CF2},
Corollary 3.4, $\rho$ is injective for \emph{all} non-hyperelliptic curves $X$.
Our result shows that if $X$ is general, then the image of $\rho$ is 
transversal to the kernel of the multiplication map ${\rm Sym}^2(H^0(K_X^{\otimes 2}))\to
H^0(K_X^{\otimes 4})$.

The proof of Theorem \ref {thm1} is by degeneration to a reducible nodal curve for which the limit of $\mu$,
described in \S\ref{2ndGauss}, has maximal rank. The theorem then follows by upper semicontinuity. 
We do not use graph curves here, i.e.\ the curves exploited in \cite{CHM},
because for them the limit of $\mu$ is more difficult to understand.
We used instead a general \emph{binary curve}, i.e.\ a stable curve of genus $g$
consisting of two rational components meeting at $g+1$ points,
which are general on both components.
For such a curve $C$ we explicitly write down the ideal $\I_2(K_C)$ in \S\ref{binary}.
In \S\ref{nontorsion} we describe the 2nd Gaussian map for $C$ modulo torsion,
and then, in \S\ref{Storsion}, we deal with the torsion part.
By direct computations performed with Maple (the script is presented
and commented in the Appendix), we verified the injectivity for a general binary curve
of genus $g \le 17$ and the surjectivity for $g=18$.
Finally, in \S\ref{Sinduction}, we proceed by induction on $g$
to complete the argument for $g \ge 19$.

The behaviour of $\mu$, and its connection
with the curvature of $\M_g$ in $\A_g$,
indicates possible relations of the surjectivity of $\mu$
with the Kodaira dimension of $\M_g$ being non-negative.
This, we think, would be a great subject for future research.
Also interesting is the study of the Gaussian maps
$\mu_k, w_k$ for higher values of $k$. The maps $\mu_k$
are related to higher fundamental forms of the Torelli immersion of $\M_g$ in $\mathcal A_g$
at a non-hyperelliptic point.  Are these maps also of maximal rank for a general curve? 

\medskip
In this paper we work over the complex field and we will use standard notation
in algebraic geometry. In particular, if $X$ is a Gorenstein curve,
$\Omega_X^1$ will denote its sheaf of K\"ahler differentials
and $K_X$ will denote its dualizing sheaf or canonical bundle, or a canonical divisor.
In general, we will indifferently use  sheaf, bundle or divisor notation.
We will often write $H^i(\cL)$ instead of $H^i(X,\cL)$ for cohomology spaces. 

\medskip

The second author wishes to thank G.~P.~Pirola for having mentioned to him
the problem solved in this paper and both G.~P.~Pirola and P.~Frediani
for discussions on this subject.

\section{The 2nd Gaussian map for a stable curve}\label{2ndGauss}

Let $X$ be a stable curve of genus $g$.
We will denote by $\I_2(K_X)$ the vector space of forms of degree 2
vanishing on the canonical model of $X$.
If $X$ is smooth,
the 2nd Gaussian map  $\mu\colon \I_2(K_X) \to H^0(X, K_X^{\otimes4})$
is locally defined as follows.

Fix a basis $\{ \omega_i \}$ of $H^0(K_X)$,
and write it in a local coordinate $z$  as 
$\omega_i = f_i(z) \,d z$.
Let $Q \in \I_2(K_X)$, with $Q = \sum_{i,j} s_{ij} \omega_i\otimes \omega_j$,
 the matrix $(s_{ij})$ being symmetric.
Since $\sum_{i,j} s_{ij} f_i f_j \equiv 0$,
one has  $\sum_{i,j} s_{ij} f'_i f_j \equiv 0$.
The local expression of $\mu(Q)$ is then (cf., e.g., \cite{CF2})

\begin{equation}\label{mu2X}
\mu(Q) = \sum_{i,j} s_{ij} f''_i f_j \,(d z)^4
= - \sum_{i,j} s_{ij} f'_i f'_j \,(d z)^4.
\end{equation}

If $X$ is nodal, one can similarly define
the 2nd Gaussian map
\[
\mu\colon \I_2(K_X) \to H^0(X, \Sym^2(\Omega^1_X) \otimes K_X^{\otimes2})
\]
which is locally defined in a similar way as in \eqref{mu2X}. Precisely,
let $\{ \omega_i \}$ be a basis of $H^0(K_X)$.
In local coordinates, we can write $\omega_i = f_i \,\xi$,
where $f_i$ is a regular function and $\xi$ is a local generator of the canonical bundle $K_X$.
Then $\mu$ is locally defined by

\begin{equation}\label{mu2a}
\mu(Q) = - \sum_{i,j} s_{ij} \,d f_i \,d f_j \,\xi^{\otimes2}.
\end{equation}

Given a flat degeneration over a disc of a general curve to a stable curve $X$,
the 2nd Gaussian map for $X$ is the flat limit of the 2nd Gaussian map for the general curve.

It is useful to describe in some detail the space $H^0(X, \Sym^2(\Omega^1_X) \otimes K_X^{\otimes2})$.
First remark that $\Sym^2(\Omega^1_X)$ has torsion $T$ supported at the nodes of $X$.  So we have
a short exact sequence
\[
0 \to T \to \Sym^2(\Omega^1_X) \to \F_X \to 0,
\]
where $\F_X$ is a non-locally free, rank 1, torsion free sheaf on $X$.

\begin{lemma}\label{torsionatnode}
(a) For every node $p$ of $X$, $T_{p}$ is a 3-dimensional vector space;
if the local equation of $X$ around $p$ is $xy=0$,
then $T_{p}$ is
spanned by $d x\,d y$, $x\,d x\,d y$ and $y\,d x\,d y$.

\noindent
(b) If $X_i$ are the irreducible components of the normalization $\pi\colon\tilde X\to X$ of $X$, one has
\[
\F_X \cong \sideset{}{_{\!i}}\bigoplus \pi_* K_{X_i}^{\otimes 2}.
\]
\end{lemma}

\begin{proof}
Since $y \,d x = -x \,d y$,
a local section of $\Sym^2(\Omega^1_X)$ around a node $xy=0$ can be uniquely written as
$f(x) \,(d x)^2 + g(x,y) \,d x\,d y + h(y) \,(d y)^2$, where $g(x,y)$ is linear.
Then (a) is a local computation and (b) follows from (a).
\end{proof}

As a consequence, since $K_{X|X_i}=K_{X_i}(D_i)$
where $D_i$ be the divisor of nodes on $X_i$,
one has
\begin{equation}\label{Sym2}
H^0(X, \Sym^2(\Omega^1_X) \otimes K_X^{\otimes2}) \cong T \oplus
\sideset{}{_{\!i}}\bigoplus H^0(X_i, K_{X_i}^{\otimes 4}(2D_i)).
\end{equation}
where $T\cong \C^{3\delta}$, with $\delta$ the number of nodes of $X$.

\section{Canonical binary curves}\label{binary}

Let $[x_1,\ldots,x_g]$ be homogenous coordinates in $\PP^{g-1}$, $g\ge 3$.
Let $p_h=[0,\ldots,0,1,0,\ldots,0]$, with 1
at the $h$-th place, $1 \le h \le g$, 
be the coordinate points and $u=[1,1,\ldots,1]$ the unit point.
Take $C_1,C_2$ two distinct rational normal curves in $\PP^{g-1}$
passing through $p_h$, $1 \le h \le g$,  and $u$. Then
$C_1, C_2$ intersect transversally at these $g+1$ points and have no
further intersection. 

We may and will assume that $C_k$, $k=1,2$, is the closure of the image
of the map $f_k$ given by
\begin{align} \label{deff(t)}
t  &\mapsto
f_k(t)=\left[ \frac{1}{t-\alpha_{k,1}}\,, \frac{1}{t-\alpha_{k,2}}\,, \ldots,
\frac{1}{t-\alpha_{k,g}} \right],
\end{align}
where $\alpha_{k, i} \in \C$, $k=1,2$, $i=1,\ldots,g$.
In particular, $f_k(\alpha_{k,h})=p_h$, $h=1,\ldots,g$, and $f_k(\infty)=u$.
For our purposes, the $\alpha_{k,i}$'s will be general in $\C$.
Actually, we will often consider them as indeterminates on $\C$. 

The curve $C=C_1\cup C_2$ is the limit of a general canonical curve
$X\subset\PP^{g-1}$ of genus $g$,
and $C$ is canonical too, i.e.\ $\mathcal O_C(1)\cong K_C$.
The curve $C$ is usually called a \emph{canonical binary curve}. 

\begin{proposition}\label{projnormal}
A canonical binary curve $C=C_1\cup C_2$ is projectively normal.
\end{proposition}

\begin{proof} The assertion is trivial for $g=3$, which is 
the minimum allowed value of $g$. So we may assume $g\geq 4$.
By Theorem 1.2 in \cite{S}, it suffices to show that
there are $g-2$ smooth points of $C$ spanning a $\PP^{g-3}$
which meets $C$ scheme-theoretically at these $g-2$ points only.
Choose $g-2$ general points on $C_1$
and let $\Lambda\cong\PP^{g-3}$ be their span.
This meets transversally $C_1$ at these points.
We claim that $\Lambda$ does not meet $C_2$.
Otherwise choose $g-4$ general points on $C_1$ and project $C$ down  to $\PP^3$
from their span. The image of $C_1$ is a rational normal cubic
$\Gamma_1$, whereas $C_2$ projects birationally (cf.\ \cite{CC})
to a non-degenerate rational curve $\Gamma_2$ 
of degree larger than 3, thus $\Gamma_1$
and $\Gamma_2$ are distinct. Moreover the general secant line
to $\Gamma_1$ would meet $\Gamma_2$, which is impossible by
the trisecant lemma (see the \emph{focal proof} in \cite {ChC}).
\end{proof}

\begin{remark}
The only information that we will need from the above proposition
is that $C$ is quadratically normal, which is equivalent to
\[
\dim (\I_2(K_C)) = \binom{g-2}{2}.
\]
The simple argument in the proof of Proposition \ref{projnormal}
relies on Schreyer's result, which requires a careful analysis,
following the classical approach of Petri. The same result 
would follow by proving that the general hyperplane section
of $C$ verifies the general position theorem (see \cite {ACGH}, p.\ 109).
This may be proved with the same argument as above, but we do not dwell on this here.

In case $C$ is a general binary curve, it is quite simple to prove that $C$ is quadratically normal.
One way is to remark that the general trigonal binary curve is quadratically normal.
For example, if $g=2h$, embed $\FF_0$ in $\PP^{g-1}$ via the linear system
of curves of type $(1,h-1)$. The general trigonal binary curve is the union of the images
of a general curve of type $(1,h)$ and of a general curve of type $(2,1)$.
The case $g$ odd is similar and is left to the reader.
\end{remark}

We are now  interested  in explicitly describing  the vector space $\mathcal I_2(K_C)$ of degree two forms vanishing on $C$,
i.e.\ the domain of the map $\mu$ for $C$. The analysis we are going to make will provide another proof that the general binary curve $C$ is quadratically normal.

For $k=1,2$, set
\begin{align}\label{Akt}
A_k(t) &=
        \prod_{i=1}^g( t-\alpha_{k, i}).
\end{align}
For each $h=0,\ldots,g$, the coefficients $c_{k,h}$ of $t^{g-h}$ in $A_k(t)$ are,
up to sign, the elementary symmetric functions
\begin{align*}
c_{k,0} &= 1,
  &
c_{k,h} &= (-1)^h \!\!\sum_{1 \le i_1 < i_2 < \cdots < i_h \le g}
   \!\!\alpha_{k,i_1}\alpha_{k,i_2}\cdots\alpha_{k,i_h}.
\end{align*}
Note that
the index $h$ is the degree of $c_{k,h}$ as a polynomial in 
the $\alpha_{k,i}$'s .

Fix $k\in\{1,2\}$.
Since  $C_k$ passes through the coordinate points, the equation of a quadric $Q\subset\PP^{g-1}$ containing $C_k$ has the form
\begin{equation}\label{eqQ}
\sum_{1\le i<j \le g} s_{ij}x_ix_j = 0,
\end{equation}
with the conditions
\begin{align*}
P_k(t):=\sum_{1\le i<j \le g} \frac{A_k(t)}{(t-\alpha_{k,i})(t-\alpha_{k,j})}\,s_{ij}
= \sum_{n=0}^{g-2} P_{k,n} t^{n} \equiv 0,
\end{align*}
where $P_k(t)$ is a polynomial in $t$ of degree $g-2$
whose coefficients are linear polynomials $P_{k,n}(s_{ij})$ in the $s_{ij}$'s,
$n=0,\ldots,g-2$.
By expanding the product $A_k(t)$ 
one sees that the coefficients $p_{k,h;i,j}$
of $s_{ij}$ in $P_{k,g-2-h}$, $h=0,\ldots,g-2$, are
\begin{align}\label{phij}
   p_{k,0;i,j} &= 1,
&  p_{k,1;i,j} &= -\!\!\sum_{i_1 \ne i,j} \alpha_{k,i_1},
&  p_{k,h;i,j} &= (-1)^{h} \!\!\!\!\!\sum_{\substack{i_1 < i_2 < \cdots < i_h\\ \text{all} \ne i,j}}
   \!\!\!\!\!\alpha_{k,i_1}\alpha_{k,i_2}\cdots\alpha_{k,i_h},
\ 2 \le h \le g-2,
\end{align}
namely the elementary symmetric functions, removing the $i$ and $j$ terms, up to sign.
Again the index $h$ coincides with the degree of $p_{k,h;i,j}$
as a homogeneous polynomial in the $\alpha_{k,i}$'s.

Consider also the polynomials
\begin{equation*}
Q_{k,n}(s_{ij}):=
\sum_{1\le i<j \le g}
\Biggr( \sum_{m=0}^{g-2-n} \alpha_{k,i}^{m} \alpha_{k,j}^{g-2-n-m} \Biggr)
  s_{ij},
\qquad n=0,\ldots,g-2,
\end{equation*}
and let $q_{k,h;i,j}=\sum_{m=0}^{h} \alpha_{k,i}^{m} \alpha_{k,j}^{h-m}$
be the coefficient of $s_{ij}$ in $Q_{k,g-2-h}$, $h=0,\ldots,g-2$.
Also in this case the index $h$ coincides with the degree of $q_{k,h;i,j}$
as a homogeneous polynomial in the $\alpha_{k,i}$'s.

\begin{remark}\label{q_h}
The coefficient $q_{k,h;i,j}$ of $s_{ij}$ in $Q_{k,g-2-h}$
can be recursively computed by
\begin{align*}
  q_{k,0;i,j}  &=  1,
& q_{k,1;i,j}  &=  \alpha_{k,i}+\alpha_{k,j},
& q_{k,h;i,j}  &=  q_{k,1;i,j} q_{k,h-1;i,j} - \alpha_{k,i}\alpha_{k,j} q_{k,h-2;i,j},
\quad 2 \le h \le g-2.
\end{align*}
Note that all the monomials $\alpha_{k,j}^{m}\alpha_{k,i}^{h-m}$, $m=0,\ldots,h$,
in particular $\alpha_{k,i}^h$ and $\alpha_{k,j}\alpha_{k,i}^{h-1}$,
appear in $q_{k,h;i,j}$ with coefficient $1$.
Note also the recursive formula
\begin{equation}\label{qhij2}
q_{k,h;i,j}=\alpha_i q_{k,h-1;i,j}+\alpha_j^h, \quad 1 \le h \le g-2.
\end{equation}
\end{remark}

We will need the following lemma:

\begin{lemma}\label{PQ}
Fix $k\in\{1,2\}$.
For each $n=0,\ldots,g-2$, one has
\begin{equation}\label{P-Q}
P_{k,n} = \sum_{m=0}^{g-2-n} c_{k,m} Q_{k,n+m}.
\end{equation}
In particular, the linear system
\begin{equation}\label{Phij}
P_{k,n}(s_{ij})=0,\qquad n=0,\ldots,g-2,
\end{equation}
in the $s_{ij}$'s  is equivalent to the linear system
\begin{equation}\label{Qhij}
Q_{k,n}(s_{ij})=0,
\qquad n=0,\ldots,g-2.
\end{equation}
\end{lemma}

\begin{proof}
One has $P_{k,g-2}=Q_{k,g-2}$ and $P_{k,g-3}=Q_{k,g-3}+c_{k,1} Q_{k,g-2}$.
Next we proceed by induction: formula \eqref{P-Q} is equivalent to
\begin{equation}\label{pkhij}
p_{k,h;i,j} = \sum_{l=0}^{h} c_{k,l} q_{k,h-l;i,j},
\qquad \text{for} \quad  h=0,\ldots,g-2.
\end{equation}
For $h=0,1$, \eqref{pkhij} clearly holds.
Since the index $k$ is fixed, we omit it.
For $2 \le h \le g-2$, one has
\begin{align*}
p_{h;i,j}-c_{h} q_{0;i,j}  &=
(\alpha_{i}+\alpha_{j})p_{h-1;i,j}-\alpha_{i}\alpha_{j}p_{h-2;i,j} \stackrel{\text{(by induction)}}{=} \\
&= c_{h-1}q_{1;i,j}+\sum_{l=0}^{h-2}c_{l}(q_{h-l-1;i,j}q_{1;i,j}
-\alpha_{i}\alpha_{j}q_{h-l-2;i,j}) =  \sum_{l=0}^{h-1} c_{l} q_{h-l;i,j},
\end{align*}
which proves \eqref{pkhij} and therefore \eqref{P-Q}.
Since $c_{k,0}=1$, the base change matrix between the $Q_{k,n}$'s and the $P_{k,n}$'s
is unipotent triangular, hence it is invertible. The equivalence between \eqref{Phij}
and \eqref{Qhij} follows.
\end{proof}

Next we can give the announced description of  $\mathcal I_2(K_C)$. 

\begin{proposition}\label{PQ2} Let $g\ge 3$.
For a general choice of $\alpha_{k,i}$, $1\le k\le 2$, $1\le i\le  g$, one has that

\begin{enumerate}[(a)]

\item the linear system \eqref{Qhij} has maximal rank $g-1$;

\item the linear system
\begin{equation}\label{Q12}
Q_{1,0}(s_{ij})=\cdots=Q_{1,g-2}(s_{ij})=
Q_{2,0}(s_{ij})=\cdots=Q_{2,g-3}(s_{ij})=0,
\end{equation}
has maximal rank $2g-3$.
\end{enumerate}
\end{proposition}

\begin{proof}
(a) Since the index $k$ is fixed, we drop it here. Let us consider the matrix 
\[
U:=U(\alpha_1,\ldots,\alpha_g)=(q_{h;i,j})_{0\le h\le g-2,1\le i<j\le g}
\]
of size $(g-1)\times \binom{g}{2}$, where
the pairs $(i,j)$ are lexicographically ordered. We have to prove that 
there is a minor of $U$ of order $g-1$ 
which is not identically zero. We show this for 
the minor 
$D:=D(\alpha_1,\ldots,\alpha_g)$ determined
by the first $g-1$ columns, indexed by $(1,i)$ with $2\le i\le  g$.
This is true if $g=3$, so we proceed by induction on $g$. 
Look at $D$ as a polynomial in $\alpha_g$:
it has degree $g-2$ and
the coefficient of $\alpha_g^{g-2}$ is
$D(\alpha_1,\ldots,\alpha_{g-1})$ (cf.\ Remark \ref {q_h}),
which is non-zero by induction. 
This proves the assertion.

Equivalently, 
by subtracting from each row the previous one multiplied by $\alpha_1$
and using \eqref{qhij2} (cf.\ Remark \ref {q_h}),
one sees that $D$
is the Vandermonde determinant
$V(\alpha_2,\ldots,\alpha_g)=\prod_{2\le i<j\le g}(\alpha_j-\alpha_i)$
of $\alpha_2,\ldots,\alpha_g$.

\medskip
\noindent
(b) We use the same idea of the proof of (a).
Form a matrix  $Z:=Z(\alpha_{k,i})_{1\le k\le 2, 1\le i<j\le g}$
of size $(2g-3)\times \binom{g}{2}$ by concatenating vertically $U$ (for $k=1$)
and the matrix
\[
W:=W(\alpha_{2,1},\ldots,\alpha_{2,g})=(q_{2,h;i,j})_{1\le h\le g-2,1\le i<j\le g}.
\]
It suffices to prove that the minor 
$M:=M(\alpha_{k,i})_{1\le k\le 2, 1\le i<j\le g}$ of $Z$ determined
by the first $2g-3$ columns, indexed by $(1,i), (2,j)$ with $2\le i\le  g$ and
$3\le j\le g$, is 
not identically zero as a polynomial in the $\alpha_{k,i}$'s. 
This is clearly true for $g=3$, so we proceed by induction on $g$. 
Look at $M$ as a polynomial in $\alpha_{1,g}$ and $\alpha_{2,g}$:
one sees that the monomial $\alpha_{1,g}^{g-2}\alpha_{2,g}^{g-3}$
appears in $M$ with the coefficient
$(\alpha_{2,2}-\alpha_{2,1})M(\alpha_{k,i})_{1\le k\le 2, 1\le i<j\le g-1}$,
which is non-zero by induction,
proving the assertion.

Equivalently, looking at $M$ as a polynomial in $\alpha_{1,1}$,
one sees that the coefficient of the monomial $\alpha_{1,1}^{g-2}$
is the product of the two Vandermonde determinants
$V(\alpha_{2,2},\ldots,\alpha_{2,g}) V(\alpha_{1,3},\ldots,\alpha_{1,g})$.
\end{proof}

\begin{remark}
The solutions of the linear system \eqref{Qhij}, as well as those of \eqref{Phij},
give us the quadrics containing the rational normal curve $C_k$,
whereas the solutions of \eqref{Q12} give us the quadrics in $\I_2(K_C)$ 
for the binary curve $C=C_1 \cup C_2$.
\end{remark}

\section{Binary curves: the 2nd Gaussian map modulo torsion}\label{nontorsion}

Let $C=C_1 \cup C_2$ be a general binary curve.
In this section we will consider the composition $\nu$ of the 2nd Gaussian map
for $C$ with the projection to the non-torsion part of
$H^0(C,\Sym^2 (\Omega^1_C) \otimes K_C^{\otimes2})$ 
(cf.\ formula \eqref{Sym2} in \S \ref {2ndGauss}).
Specifically, for $k=1,2$, we will look at the map
\[
\nu_k\colon \I_2(K_C) \to H^0(C_k, K_{C_k}^{\otimes4}(2D_k))
\]
where $D_k$ is a divisor of degree $g+1$ on $C_k$,
therefore $\nu=\nu_1 \oplus \nu_2$ and
\[
H^0(C_k, K_{C_k}^{\otimes4}(2D_k)) \cong H^0(\PP^1, \cO_{\PP^1}(2g-6)).
\]

The map $\nu_k$ can be explicitly written down, by taking into account \eqref{mu2a}
and the description of the ideal $\I_2(K_C)$ (see \S \ref {binary}).
Precisely, the let $Q\in \I_2(K_C)$ be of the form \eqref{eqQ}
where the $s_{ij}$'s are solutions of \eqref{Q12}.
Then
\[
\nu_k(Q)=\sum_{1\le i< j \le g} \frac{1}{(t-\alpha_{k,i})^2(t-\alpha_{k,j})^2}\, s_{ij} (d t)^4
\in H^0(C_k, K_{C_k}^{\otimes4}(2D_k)).
\]
To look at this as a section of $H^0(\PP^1, \cO_{\PP^1}(2g-6))$,
we multiply by $A_k^2(t)$. Hence
\begin{equation}\label{Rk}
\nu_k(Q)=\sum_{1\le i< j \le g} \frac{A_k^2(t)}{(t-\alpha_{k,i})^2(t-\alpha_{k,j})^2}\, s_{ij}=:R_k(t)
\end{equation}
is a polynomial in $t$ whose apparent degree is $2g-4$,
but its coefficient of degree $2g-4$
is $P_{k,g-2}$ and the one of degree $2g-5$ is proportional to $P_{k,g-3}$,
hence they vanish and $R_k(t)$ has actual degree $2g-6$.

Using this explicit description \eqref{Rk} of $\nu$, we asked Maple to compute its rank
for low values of $g$ (see the Appendix for Maple scripts). The result is the following:

\begin{proposition}\label{nu11}
The map $\nu$ has maximal rank for $g \le 18$,
namely $\nu$ is injective for $g \le 10$ and it is surjective for $11 \le g \le 18$.
\qed
\end{proposition}

\begin{corollary}
The 2nd Gaussian map $\mu$ is injective for the general curve of genus $g \le 10$.
\qed
\end{corollary}

\section{Binary curves: the torsion}\label{Storsion}

Let $C=C_1\cup C_2$ be a general binary curve as in \S\ref{binary}.
In \eqref{deff(t)} we may replace $f_k$, $1\le k\le 2$, with
\begin{equation}\label{phikit}
A_k(t) f_k(t)=\left[
\phi_{k,1}(t),\ldots,\phi_{k,g}(t)
\right],
\qquad  \phi_{k,i}(t)=\frac{A_k(t)}{(t-\alpha_{k,i})}.
\end{equation}

Now we consider the restriction $\tau$ of the 2nd Gaussian map
for $C$ to ${\rm ker}(\nu)$, which lands in the torsion part $T$ of
$H^0(C,\Sym^2 (\Omega^1_C) \otimes K_C^{\otimes2})$,
cf.\ formula \eqref{Sym2}.
By taking into account Lemma \ref{torsionatnode}, (a),
a direct computation shows that
the composition of $\tau$ with the projection
on the torsion part $T_{p_h}$
at the coordinate point $p_h$ is as follows:
if $Q\in {\rm ker}(\nu)$ is of the form \eqref{eqQ},
then $Q$ is mapped to
{\small\begin{equation}\label{torsionatph}
dx\,dy\sum_{i\ne j} s_{ij} \phi'_{1,i}(\alpha_{1,h}) \phi'_{2,j}(\alpha_{2,h})
+2x\,dx\,dy \sum_{i\ne j} s_{ij} \phi''_{1,i}(\alpha_{1,h}) \phi'_{2,j}(\alpha_{2,h})
+2y\,dx\,dy \sum_{i\ne j} s_{ij} \phi'_{1,i}(\alpha_{1,h}) \phi''_{2,j}(\alpha_{2,h}),
\end{equation}}%
where $s_{ji}=s_{ij}$ and $x,y$ are local coordinates around $p_h$ such that $C_1\colon y=0$
and $C_2\colon x=0$.
The description of the torsion at the unitary point $u$ is similar.
Replace $f_k$ by the parametrization
$
\frac{1}{t}\,f_k(\frac{1}{t}).
$
Again a direct computation shows that 
the composition of $\tau$ with the projection on $T_{u}$ is
\begin{equation}\label{torsionatu}
Q \mapsto
dx\,dy\sum_{i\ne j} s_{ij} \alpha_{1,i}\alpha_{2,j}
+2x\,dx\,dy \sum_{i\ne j} s_{ij} \alpha_{1,i}^2\alpha_{2,j}
+2y\,dx\,dy \sum_{i\ne j} s_{ij} \alpha_{1,i}\alpha_{2,j}^2,
\end{equation}
where $s_{ji}=s_{ij}$ and $x,y$ are local coordinates around $u$ such that $C_1\colon y=0$
and $C_2\colon x=0$.

Consider the following commutative diagram with exact rows
\begin{equation}\label{tau}
\raisebox{4ex}{\xymatrix@C-1ex@R-1ex{%
0 \ar[r] & T \ar[r]
  & H^0(C,\Sym^2(\Omega^1_C) \otimes K_C^{\otimes2}) \ar@{->>}[r]
  & H^0(C_1,K_{C_1}^{\otimes2}(2)) \oplus H^0(C_2,K_{C_1}^{\otimes2}(2)) \cong H^0(\F_C) \\
0 \ar[r] & \ker(\nu) \ar[u]^-{\tau} \ar[r] & \I_2(K_C) \ar[u]^-{\mu} \ar[ur]!UR_-{\nu}
}}
\end{equation}

We asked Maple to compute the rank of the map $\tau$ for $11 \le g \le 18$
(see the script in the Appendix).
Taking into account diagram \eqref{tau},
the result is the following:

\begin{proposition}\label{prop<18}
Let $C$ be a general binary curve of genus $g$.
The maps $\tau$ and $\mu$ have maximal rank for $g \le 18$,
namely they are injective for $g \le 17$ and surjective for $g=18$.
\qed
\end{proposition}

\begin{corollary}\label{corg18}
The map $\mu$ is injective for the general curve
of genus $g \le 17$, and surjective for $g=18$.
\qed
\end{corollary}

\section{The induction step}\label{Sinduction}

In this section we prove the main result of this paper,
namely the surjectivity of the 2nd Gaussian map $\mu$
for the general curve of genus $\ge 18$.

Let $C\subset \PP^{g-1}$ be a nodal canonical curve
and let $p\in C$ be a node.
Let $\tilde C \to C$ be the partial normalization of $C$ at $p$,
and let $p_1,p_2 \in \tilde C$ be the points over $p$.
Note that the projection from $p$ maps $C$ to the canonical model of $\tilde C$ in $\PP^{g-2}$
and we will assume that this induces an isomorphism of $\tilde C$ to its canonical model.
Consider the following diagrams 
\begin{equation}\label{chi}
\begin{aligned}
&
\xymatrix@C-1ex@R-1ex{%
 0\ar[r] & H^0(\F_{\tilde C}) \ar@{^{(}->}[r] & H^0(\F_{C}) \ar@{->>}[r] & \cO_{2p_1}\oplus\cO_{2p_2}  \\
0 \ar[r] & \I_2(K_{\tilde C}) \ar@{^{(}->}[r] \ar[u]^{\tilde\nu} & \I_2(K_{C}) \ar[u]^{\nu} \ar[ur]!UR_-{\chi}
}
&&
\xymatrix@C-1ex@R-1ex{%
 0\ar[r] & \tilde T \ar@{^{(}->}[r] & T \ar@{->>}[r] & T_p  \\
0 \ar[r] & \ker(\tilde\nu) \ar@{^{(}->}[r] \ar[u]^{\tilde\tau} & \ker(\nu) \ar[u]^{\tau} \ar[ur]_-{\tau_p}
}
\end{aligned}
\end{equation}
where $\tilde T$ is the torsion subsheaf of ${\rm Sym}^ 2(\Omega^ 1_{\tilde C})$, 
 $\nu, \tau$ are the maps of diagram \eqref{tau} for the curve $C$ and $\tilde\nu, \tilde\tau$ are the corresponding ones for $\tilde C$.
Diagrams \eqref{chi} are commutative and the horizontal sequences are exact, hence the next lemma is clear:

\begin{lemma}
If $\tilde\nu$ and $\chi$ [$\tilde\tau$ and $\tau_p$, resp.] are surjective, then $\nu$ [$\tau$, resp.] is also surjective.
\qed
\end{lemma}

We apply this to prove:

\begin{theorem}\label{thm>18}
If $C=C_1\cup C_2$ is a general binary curve of genus $g \ge 18$,
then $\mu$ is surjective for $C$.
\end{theorem}

\begin{proof}
The case $g=18$ has been done
by a direct computation, cf.\  Proposition \ref{prop<18}.
We then proceed by induction on $g$:
the commutativity of the diagram \eqref{tau} and the previous lemma
show that it is enough to prove the surjectivity of $\chi$ and $\tau_p$,
where $p$ is a node of $C$. We will do this for $p=u$ the unitary point.

In this situation,
the map $\nu$ is the one $\nu_1\oplus\nu_2$ considered in \S\ref{nontorsion}.
Therefore $\chi=\chi_1\oplus\chi_2$,
where $\chi_k$ is the composition of $\nu_k$ with the restriction to $\cO_{2p_k}$, $k=1,2$.
In local coordinates, $\chi_k(Q)$ is the pair formed by the constant term and the coefficient of
the degree-one term of the Taylor expansion around $p$ of the polynomial $\nu_k(Q)$.
In \S\ref{nontorsion} we computed $\nu_k$
using a local coordinate $t$ on $C_k$.
In this coordinate, the point $p=[1,\ldots,1]$ corresponds to $t=\infty$.
Therefore, if $Q\in \I_2(K_C)$ is of the form \eqref{eqQ}, with  the $s_{ij}$'s satisfying \eqref{Q12},
then $\chi_k(Q)$ is the pair of coefficients of the highest degrees 
$2g-6$ and $2g-7$ of the polynomial $\nu_k(Q)$, i.e.\ of the polynomial $R_k(t)$ given in \eqref{Rk}.
We denote by $R_{k,2g-6}$ and $R_{k,2g-7}$
these coefficients, which are linear polynomials in the $s_{ij}$'s.
We will now compute them.

Fix the index $k$ and omit it.
By expanding $A^2$ in \eqref{Rk}, one sees that
the coefficient of $s_{ij}$ in $R_{2g-6}$ is
\[
4p_{2;i,j}+\sum_{i_1\ne i,j} \alpha_{i_1}^2
=4p_{2;i,j}+n_{2}-(\alpha_{i}^2+\alpha_{j}^2),
\]
where $n_{2}=\sum_{m=1}^g \alpha_{m}^2$ is independent of $i,j$,
and $p_{2;i,j}$ is the coefficient of $s_{ij}$ in $P_{k,g-4}$,
cf.\ \eqref{phij}.
By \eqref{Phij}, this means that
\[
R_{2g-6} =4P_{g-4}+n_2 P_{g-2}-\sum_{i<j}(\alpha_{i}^2+\alpha_{j}^2)s_{ij}
=-\sum_{i<j}(\alpha_{i}^2+\alpha_{j}^2)s_{ij}.
\]
Similarly, one sees that the coefficient of $s_{i,j}$ in $R_{2g-7}$ is twice
\[
  -4p_{3;i,j}-\!\!\sum_{\substack{i_1\ne i_2\\\text{both}\ne i,j}}
           \!\! \alpha_{i_1}^2\alpha_{i_2}
=-4p_{3;i,j}+n_3+n_2q_{1;i,j}-c_1(\alpha_i^2+\alpha_j^2)
-(\alpha_i^3+\alpha_j^3)-q_{3;i,j}.
\]
where $n_{3}=-\sum_{m=1}^g \alpha_{m}^3$ is independent of $i,j$.
Therefore, taking into account \eqref{Phij} and \eqref{Qhij}, one has
\[
R_{2g-7}=-2c_1R_{2g-6}-2\sum_{i<j}(\alpha_i^3+\alpha_j^3)s_{ij}.
\]

Form the matrix $Y:=Y(\alpha_{k,i})_{1\le k\le 2, 1\le i<j\le g}$ of size $(2g+1)\times \binom{g}{2}$ obtained by concatenating vertically the matrix
$Z$ in the proof of Proposition \ref {PQ2}, (b), and the matrix of size $4\times \binom{g}{2}$  whose rows are
$(\alpha_{k,i}^h+\alpha_{k,j}^h)_{1\le i<j\le g}$, with $1\le k\le 2$, $2\le h\le 3$. In order 
to accomplish the proof that $\chi$ is surjective we have to prove that there is a minor of order $2g+1$ of $Y$ which is not identically zero.
We will do this for the minor $N:=N(\alpha_{k,i})_{1\le k\le 2, 1\le i<j\le g}$
determined by the first $2g+1$ columns, indexed by $(1,i)$, $(2,j)$, $(3,\ell)$, with
$2\le i\le g$, $3\le j\le g$, $4\le \ell \le 7$. This is non-zero for $g=7$: we verified
this with Maple (see the script in the Appendix). Then we proceed by induction on $g$ and
we assume $g\ge 8$. 
The argument here is the same as the one in the proof of Proposition \ref {PQ2}, (b).
Look at $N$ as a polynomial in $\alpha_{1,g}$ and $\alpha_{2,g}$:
the monomial $\alpha_{1,g}^{g-2}\alpha_{2,g}^{g-3}$ appears in $N$
with coefficient
$(\alpha_{2,2}-\alpha_{2,1})N(\alpha_{k,i})_{1\le k\le 2, 1\le i<j\le g-1}$,
which is non-zero by induction, proving that $\chi$ is surjective.

It remains to show that $\tau_p$ is surjective. 
This could be seen with a quick monodromy argument, 
on which however we do not dwell, preferring to  present instead an argument in 
the same style as the  ones we made so far. 

Recall that $\ker(\nu)$ is defined in $ \I_2(K_C)$ by the vanishing of the polynomials $R_k(t)$, $k=1,2$, whose coefficients of degree at most $2g-8$ are polynomials in the $\alpha_{k,i}$'s of degree at least 4.
By the description of the torsion at the unitary point  given in \eqref{torsionatu}, we need to show the rank maximality of the matrix $Y'=Y'(\alpha_{k,i})_{1\le k\le 2, 1\le i<j\le g}$  of size $(2g+4)\times \binom{g}{2}$ obtained by concatenating vertically the above matrix $Y$ and the matrix of size $3\times \binom{g}{2}$  whose rows are $(\alpha_{1,i}\alpha_{2,j}+\alpha_{1,j}\alpha_{2,i})_{1\le i<j\le g}$, $(\alpha_{1,i}^2\alpha_{2,j}+\alpha_{1,j}^2\alpha_{2,i})_{1\le i<j\le g}$, and $(\alpha_{1,i}\alpha_{2,j}^2+\alpha_{1,j}\alpha_{2,i}^2)_{1\le i<j\le g}$.
We claim that the minor $N'=N'(\alpha_{k,i})_{1\le k\le 2, 1\le i<j\le g}$ of $Y'$ determined by the first $2g+4$ columns, indexed by $(1,i)$, $(2,j)$, $(3,\ell)$, with
$2\le i\le g$, $3\le j\le g$, $4\le \ell \le 10$ is non-zero for $g\ge 10$. We verified the case $g=10$ with Maple (see the script in the Appendix) and the induction is the same as before because the monomial $\alpha_{1,g}^{g-2}\alpha_{2,g}^{g-3}$ appears in $N'$
again with coefficient $(\alpha_{2,2}-\alpha_{2,1})N'(\alpha_{k,i})_{1\le k\le 2, 1\le i<j\le g-1}$.
This concludes the proof that $\tau_p$ is surjective,
hence the proof of the theorem.
\end{proof}

\begin{corollary}
The 2nd Gaussian map $\mu$ is surjective for the general curve of genus $g \ge 18$.
\qed
\end{corollary}

\appendix

\section*{Appendix: Maple scripts for computations}   

We list here the Maple script we run. We will explain it afterwards:
for this purpose, we added line numbers at each five lines.

\bigskip\noindent
\verb!alpha[1]:=[3,12,21,29,37,41,43,46,54,62,65,72,81,85,89,94,97,105]:!
\\\verb!alpha[2]:=[6,18,24,36,39,42,45,52,60,63,71,80,84,86,91,96,104,108]:!
\\\verb!for g from 4 to 18 do!
\\\verb!listsij:=[seq(seq(s[i,j],j=i+1..g),i=1..g)]:!
\\\texttt{\hspace*{-4ex}5\settowidth{\spazio}{5}\hspace{-\spazio}\hspace{4ex}}%
\verb!for k from 1 to 2 do!
\\\verb!  A[k]:=mul(t-alpha[k][i],i=1..g):!
\\\verb!  R[k]:=add(add(s[i,j]*(A[k]^2)/((t-alpha[k][i])^2*(t-alpha[k][j])^2),!
\\\verb!                j=i+1..g),i=1..g):!
\\\verb!end do:!
\\\texttt{\hspace*{-5ex}10\settowidth{\spazio}{10}\hspace{-\spazio}\hspace{5ex}}%
\verb!Z:=linalg[matrix]([seq([seq(seq(add(alpha[1][i]^m*alpha[1][j]^(h-m),m=0..h),!
\\\verb!                   j=i+1..g),i=1..g)],h=0..g-2),!
\\\verb!                   seq([seq(seq(add(alpha[2][i]^m*alpha[2][j]^(h-m),m=0..h),!
\\\verb!                   j=i+1..g),i=1..g)],h=1..g-2)]):!
\\\verb!Zref:=Gausselim(Z,'r0') mod 109:!
\\\texttt{\hspace*{-5ex}15\settowidth{\spazio}{15}\hspace{-\spazio}\hspace{5ex}}%
\verb!printf("For g=%2d, one has dim I2(K)=%3d, ",g,nops(listsij)-r0):!
\\\verb!EqsKerNu:=[seq(seq(primpart(coeff(R[k],t,n)),n=0..2*g-6),k=1..2)]:!
\\\verb!K:=Gausselim(linalg[stackmatrix](Zref,!
\\\verb!             linalg[genmatrix](EqsKerNu,listsij)),'r1') mod 109:!
\\\verb!printf("dim Ker(nu)=%2d, corank(nu)=%d, ",nops(listsij)-r1,4*g-10-r1+r0):!
\\\texttt{\hspace*{-5ex}20\settowidth{\spazio}{20}\hspace{-\spazio}\hspace{5ex}}%
\verb!for k from 1 to 2 do for i from 1 to g do!
\\\verb!  phi1[k,i]:=diff(A[k]/(t-alpha[k][i]),t): phi2[k,i]:=diff(phi1[k,i],t):!
\\\verb!  for h from 1 to g do!
\\\verb!    phi1e[k,i,h]:=eval(phi1[k,i],t=alpha[k][h]):!
\\\verb!    phi2e[k,i,h]:=eval(phi2[k,i],t=alpha[k][h]):!
\\\texttt{\hspace*{-5ex}25\settowidth{\spazio}{25}\hspace{-\spazio}\hspace{5ex}}%
\verb!end do: end do: end do:!
\\\verb!for h from 1 to g do!
\\\verb!  tors[h,1]:=add(add(s[i,j]*(phi1e[1,i,h]*phi1e[2,j,h]!
\\\verb!                            +phi1e[1,j,h]*phi1e[2,i,h]),j=i+1..g),i=1..g):!
\\\verb!  tors[h,2]:=add(add(s[i,j]*(phi2e[1,i,h]*phi1e[2,j,h]!
\\\texttt{\hspace*{-5ex}30\settowidth{\spazio}{30}\hspace{-\spazio}\hspace{5ex}}%
\verb!                            +phi2e[1,j,h]*phi1e[2,i,h]),j=i+1..g),i=1..g):!
\\\verb!  tors[h,3]:=add(add(s[i,j]*(phi1e[1,i,h]*phi2e[2,j,h]!
\\\verb!                            +phi1e[1,j,h]*phi2e[2,i,h]),j=i+1..g),i=1..g):!
\\\verb!end do:!
\\\verb!tors[0,1]:=add(add(s[i,j]*(alpha[1][i]*alpha[2][j]!
\\\texttt{\hspace*{-5ex}35\settowidth{\spazio}{35}\hspace{-\spazio}\hspace{5ex}}%
\verb!                          +alpha[1][j]*alpha[2][i]),j=i+1..g),i=1..g):!
\\\verb!tors[0,2]:=add(add(s[i,j]*(alpha[1][i]^2*alpha[2][j]!
\\\verb!                          +alpha[1][j]^2*alpha[2][i]),j=i+1..g),i=1..g):!
\\\verb!tors[0,3]:=add(add(s[i,j]*(alpha[1][i]*alpha[2][j]^2!
\\\verb!                          +alpha[1][j]*alpha[2][i]^2),j=i+1..g),i=1..g):!
\\\texttt{\hspace*{-5ex}40\settowidth{\spazio}{40}\hspace{-\spazio}\hspace{5ex}}%
\verb!EqsKerTau:=[seq(seq(primpart(tors[h,l]),l=1..3),h=0..g)]:!
\\\verb!Gausselim(linalg[stackmatrix](K,linalg[genmatrix](EqsKerTau,listsij)),'r2') mod 109:!
\\\verb!printf("dim ker(tau)=%d, corank(tau)=%2d\n",nops(listsij)-r2,3*g+3-r2+r1):!
\\\verb!if g=7 then!
\\\verb!  N:=linalg[det](linalg[stackmatrix](linalg[delcols](Z,16..21),!
\\\texttt{\hspace*{-5ex}45\settowidth{\spazio}{40}\hspace{-\spazio}\hspace{5ex}}%
\verb!                 linalg[matrix]([seq(seq([seq(seq(alpha[k][i]^h+alpha[k][j]^h,!
\\\verb!                                 j=i+1..7),i=1..3)],h=2..3),k=1..2)]))):!
\\\verb!  printf("For g= 7, the minor N is congruent to %d (mod 5)\n",N mod 5):!
\\\verb!elif g=10 then!
\\\verb!  N2:=linalg[det](linalg[stackmatrix](linalg[delcols](Z,25..45),!
\\\texttt{\hspace*{-5ex}50\settowidth{\spazio}{40}\hspace{-\spazio}\hspace{5ex}}%
\verb!      linalg[matrix]([seq(seq([seq(seq(alpha[k][i]^h+alpha[k][j]^h,!
\\\verb!                                       j=i+1..10),i=1..3)],h=2..3),k=1..2)]),!
\\\verb!      linalg[matrix]([[seq(seq(alpha[1][i]*alpha[2][j]!
\\\verb!                              +alpha[1][j]*alpha[2][i],j=i+1..10),i=1..3)],!
\\\verb!                      [seq(seq(alpha[1][i]^2*alpha[2][j]!
\\\texttt{\hspace*{-5ex}55\settowidth{\spazio}{40}\hspace{-\spazio}\hspace{5ex}}%
\verb!                              +alpha[1][j]^2*alpha[2][i],j=i+1..10),i=1..3)],!
\\\verb!                      [seq(seq(alpha[1][i]*alpha[2][j]^2!
\\\verb!                              +alpha[1][j]*alpha[2][i]^2,j=i+1..10),i=1..3)]]))):!
\\\verb!  printf("For g=10, the minor N' is congruent to %d (mod 23)\n",N2 mod 23):!
\\\verb!end if:!
\\\verb!end do:!

\bigskip
In lines 1--2, we define the $\alpha_{k,i}$'s which will be used. We chose them randomly.
In line 3, we start the main loop which runs for $4 \le g \le 18$.
In line 4, we collect the unknowns $\{s_{i,j}\}_{1\le i < j \le g}$ in the list \texttt{listsij}:
there are $\binom{g}{2}$ of them.
In lines 6--8 we define the polynomials $A_k(t)$
and $R_k(t)$, cf.\ \eqref{Akt} and \eqref{Rk}.

In lines 10--13, we define the matrix \texttt{Z}
associated to the linear system \eqref{Q12},
whose solutions give us the quadrics in $\I_2(K_C)$,
cf.\ the proof of Proposition \ref{PQ2}.
In line 14, Maple computes the rank \texttt{r0} of \texttt{Z}
via Gaussian elimination, by calculating modulo 109 to speed up computations.
The resulting matrix in row echelon form is called \texttt{Zref}.
As expected by Proposition \ref{PQ2}, (b),
Maple finds $\texttt{r0}=2g-3$ for each $g=4,\ldots,18$.
In line 15, Maple prints out the genus $g$
and $\dim (\I_2(K_C))=\binom{g}{2}-\texttt{r0}=\binom{g-2}{2}$.

In line 16, we collect in \texttt{EqsKerNu} the list of equations
which determine $\ker(\nu)$,
cf.\ the definition \eqref{Rk} of $\nu$ in \S\ref{nontorsion}.
In lines 17--18, Maple computes the rank \texttt{r1} of the linear system
$\texttt{EqsKerNu} \cap\ker(\texttt{Zref})$,
again via Gaussian elimination modulo 109,
and the resulting row echelon matrix is called \texttt{K}.
Maple finds that $\texttt{r1}=\binom{g}{2}$ for $4 \le g \le 10$
and that $\texttt{r1}=6g-13$ for $11 \le g \le 18$.
Therefore the rank of $\nu$ is
$\texttt{r1}-\texttt{r0}=\binom{g-2}{2}$ for $4 \le g \le 10$,
and $=4g-10$ for $11 \le g \le 18$.
This proves Proposition \ref{nu11}.

In line 19, Maple prints out the dimension of $\ker(\nu)$
and the corank of $\nu$, that is $4g-10-\texttt{r1}+\texttt{r0}$.

In lines 20--25, we define the 1st derivative \texttt{phi1}
and the 2nd one \texttt{phi2} of the $\phi_{k,i}$'s, cf.\ \eqref{phikit}.
We then define their evaluations \texttt{phi1e}, \texttt{phi2e}
at the coordinate point $p_h$.
Using them, in lines 26--33 we compute the torsion at $p_h$, $h=1,\ldots,g$,
cf.\ \eqref{torsionatph}, and, in lines 34--39,
the torsion at the unit point $u$, cf.\ \eqref{torsionatu}.

In lines 40--41, we collect in \texttt{EqsKerTau} the equations which
determine $\ker(\tau)$ and Maple computes the rank \texttt{r2}
of $\texttt{EqsKerTau} \cap\ker(\text{\texttt{K}})$,
via Gaussian elimination modulo 109 as before.
Maple finds that $\texttt{r2}=\binom{g}{2}$ for $4 \le g \le 17$
and that $\texttt{r2}=152$ for $g=18$.
Therefore the rank of $\tau$ is
$\texttt{r2}-\texttt{r1}=(g^2-13g+26)/2$ for $11 \le g \le 17$,
and $=57$ for $g = 18$.
This proves Proposition \ref{prop<18}.

In line 42, Maple prints out the the dimension of $\ker(\tau)$
and the corank of $\tau$, that is $3g+3-\texttt{r2}+\texttt{r1}$.

Finally, in lines 43--59, Maple computes the minors $N$ (when $g=7$)
and $N'$ (when $g=10$), needed in the proof of Theorem \ref{thm>18},
and it prints out that $N\bmod{5} = 4$ and $N'\bmod{23} = 16$.

\vspace*{3mm}

\end{document}